\documentclass[12pt]{article}
\usepackage{amssymb}
\usepackage{amsthm}
\usepackage{amsmath}
\usepackage{graphicx}
\graphicspath{{./Figs/}}
\newtheorem{defi}{Definition}
\newtheorem{teo}{Theorem}
\newtheorem{lem}{Lemma}

\newcommand{\pa}{\partial}
\newcommand{\ol}{\overline}

\newcommand{\vp}{\varphi}
\newcommand{\ve}{\varepsilon}
\newcommand{\om}{\omega}
\newcommand{\be}{\begin{equation}}
\newcommand{\ee}{\end{equation}}

\newcommand{\ipd}{\stackrel{\normalfont\text{def}}{=}}
\newcommand{\const}{\operatorname{const}}
\newcommand{\dpsi}{\dot{\psi}_0}

\newcommand{\ts}{\tilde{\sigma}}

\usepackage[usenames]{color}
\begin{document}
\allowdisplaybreaks

\title{Collision of solitons for a non-homogenous version  of the  KdV equation }

\author{ Georgy~Omel'yanov\thanks{
Universidad de Sonora, Rosales y Encinas s/n, 83000, Hermosillo,
Sonora, M\'exico,\ omel@mat.uson.mx} }
\date{}
\maketitle
\begin{abstract}
We consider  KdV-type  equations with $C^1$ nonhomogeneous nonlinearities
and small dispersion $\ve$. The first result consists in the
conclusion that, in the leading term with respect to $\ve$, the
solitary waves in this model interact like KdV solitons.  Next it turned out that there
exists a very interesting scenario of instability in which the short-wave soliton remains stable
whereas a small long-wave part, generated by perturbations of original equation, turns to be unstable,
growing and destroying the leading term.  At the same time, such perturbation can eliminate the collision of solitons.
\end{abstract}
\emph{Key words}:  KdV-type  equation, soliton,
interaction,   weak asymptotics method
\emph{2010 Mathematics Subject Classification}: 35Q53, 35D30

\section{Introduction}
We consider a generalization of the KdV equation of the form:
\begin{equation}\label{1}
\frac{\pa u}{\pa t}+\frac{\pa g'(u)}{\pa x}+\varepsilon^2
\frac{\partial^3u}{\pa x^3}= 0, \; x \in \mathbb{R}^1, \; t > 0,
\end{equation}
where $g'(u)\ipd\pa g/\pa u\in C^1$ is a real function (for more detail see below) and  $\varepsilon<<1$ is a small parameter.
Such equations describe nonlinear wave phenomena in plasma physics. In particular,
for some specific  plasma states, the ion-acoustic or dust-acoustic  phenomena can be
described by the KdV-type equation (\ref{1}) with non-linearities
$g'(u)= \alpha u^{3/2}+\beta u^2$ or  $g'(u)= \alpha u^{2}+\beta u^3$, $\alpha$, $\beta =\const$  (\cite{Sch} - \cite{Rahm}, see also \cite{KudSi}).
To simplify the situation we restrict ourselves  by non-negative $u$. Moreover, we assume that uniformly in $u\geq 0$
\begin{equation}\label{2a}
c_1u^{1+\delta_1}\leq g'(u)\leq c_2u^{5-\delta_2},
\ee
where $c_i, \delta_i$ are positive constants. These restrictions  imply for $\ve=\const$ both the solvability of the Cauchy problem for (\ref{1}) and
the solution stability with respect to initial data (see \cite{Fam, Bona}).

For homogenous case $g'(u)=u^\kappa$, $\kappa>1$, it is easy to find explicit
 solitary wave solutions (see below). Moreover, as it is well known nowadays, the solitons interact
 elastically  in the integrable case ($\kappa=2$ and $3$). Almost the same is true for non-integrable homogenous case:
 the solitary waves interact elastically in the principal term in an asymptotic sense, whereas  the non-integrability
  implies the appearance of small radiation-type corrections (\cite{DanOm} - \cite{OmVal}).
  At the same time,  the character of
 the solitary wave collision remains unknown for arbitrary non-linearity.  The same is true for the solitary wave
 stability with respect to right-hand side perturbations. Our aim is to consider these open problems.

\section{Solitary wave solution}
First of all, we should determine what type of solitary waves will be  under consideration.
\begin{defi}\label{def1.1}
A  function
\begin{equation}\label{2}
u=A\om\big(\beta(x-Vt)/\varepsilon,A\big)
\end{equation}
is called "soliton type solitary wave" if  $\beta=\beta(A)$, $V=V(A)$,
and $\om=\om(\eta,A)$ are smooth functions uniformly in the parameter $A>0$, and $\om$ is an even function such that  $\om(0,\cdot)=1$,
$\om'(0,\cdot)=0$, and $\om''(0,\cdot)<0$, where the prime denotes the derivative with respect to $\eta$. Moreover, we assume that
\be\label{2b}
\om(0,\cdot)\to0\quad\text{as}\quad \eta\to\pm\infty
\ee
with an exponential rate, and $\om(\eta,\cdot)<1$ for $\eta\neq0$.
\end{defi}
\begin{teo} \label{teo1}
Let $g(u)\in C^2(u\geq0)\bigcap C^\infty(u>0)$ be such that
\be\label{2c}
g(u)=u^2g_1(u),
\ee
where the H\"{o}lder continuous function $g_1$ satisfies the conditions:
\be\label{2d}
g_1(0)=0,\qquad g_1(u)>0 \quad \text{and}\quad g'_1(u)>0\quad\text{for}\quad u>0.
\ee
Then the equation
(\ref{1}) has  a soliton type solitary wave solution such that
\be\label{2e}
\pa\om(\eta,A)/\pa A\to0\quad\text{as}\quad \eta\to0\quad\text{or}\quad\eta\to\pm\infty.
\ee
\end{teo}
To prove the statement it is enough to substitute the desired form  (\ref{2}) into the equation (\ref{1}), integrate,  and use the condition (\ref{2b}):
\be\label{2f}
-AV\om+g'(A\om)+A\beta^2\frac{d^2\om}{d\eta^2}=0.
\ee
Multiplying (\ref{2f}) by $A\om'$, taking into account  (\ref{2c}), and integrating again we obtain the first order equation:
\be\label{2g}
\beta^2\left(\frac{d\om}{d\eta}\right)^2=V\om^2\left(1-\frac{g_1(A\om)}{g_1(A)}\frac{2g_1(A)}{V}\right).
\ee
It is obvious now that our  choice of the functional class  implies the unambiguous definition of $V$, whereas $\beta$
is a free parameter. So we define them in the following form:
\be\label{2h}
V=2g_1(A),\quad \beta^2=V.
\ee
Then we obtain the final version of the equation for $\om$:
\begin{equation}\label{3}
\frac{d\om}{d\eta}=\pm\om\sqrt{1-\frac{g_1(A\om)}{g_1(A)}},
\ee
where the sign should be $-$ for $\eta>0$ and $+$ for $\eta<0$. To complete the proof it is enough to
analyze the implicit representation of $\om$ which corresponds to (\ref{3}):
\be\label{2i}
\eta=\int_\om^1z^{-1}\left(1-\frac{g_1(Az)}{g_1(A)}\right)^{-1/2}dz\quad\text{for}\quad\eta\geq0.
\ee
\qed
\medskip

\textbf{Example}.

Let $g'(u)=u^\kappa$. Then the solution of (\ref{3}) does not depend on $A$ and has the form:
\be\label{4}
\om(\eta)=\left\{\cosh\big((\kappa-1)\eta/2\big)\right\}^{-2/(\kappa-1)}, \quad V=2A^{\kappa-1}/(\kappa+1).
\end{equation}

\section{Two-soliton asymptotic solution}
\subsection{Main definitions}
Obviously, there is not any hope to find both the exact multi-soliton solution to (\ref{1}) and an
asymptotics in the classical sense. So, we will construct a
weak asymptotic solution.
The Weak Asymptotics Method (see e.g. \cite{DanOm} - \cite{KalMit} and references therein) takes
into account the fact that soliton-type solutions which are smooth for $\ve>0$ become non-smooth in the limit as $\ve\to0$. Thus, it is
possible to treat such solutions  as a mapping $\mathcal{C}^\infty
(0, T; \mathcal{C}^\infty (\mathbb{R}_x^1))$ for $\ve=\const>0$ and only
as $\mathcal{C} (0, T; \mathcal{D}' (\mathbb{R}_x^1))$ uniformly
in $\ve\ge0$. Accordingly, the remainder should be small in the
weak sense. The main advantage of the method is such that we can ignore the real shape of the colliding waves but look for (and find) exceptionally their main characteristics. For the equation (\ref{1}) they are the amplitudes and trajectories of the waves.

Similarly the famous Whitham method we define a weak asymptotic solution as a function
which satisfies some conservation laws, in fact two laws for the two-phase asymptotics.  Namely, they are the  equation (\ref{1}) itself and
\begin{equation}\label{7}
\frac{\pa u^2}{\pa t}-\frac{\pa }{\pa
x}\Big\{2\big(g(u)-ug'(u)\big)+3 \Big(\ve\frac{\partial u}{\pa
x}\Big)^2-\varepsilon^2 \frac{\partial^2u^2}{\pa x^2}\Big\}= 0.
\end{equation}
As it has been demonstrated in \cite{DanOm, Om}, the correct definition of two-soliton asymptotics is the following:
\begin{defi}\label{def2.1}
A sequence $u(t, x, \varepsilon )$, belonging to
$\mathcal{C}^\infty (0, T; \mathcal{C}^\infty (\mathbb{R}_x^1))$
for $\varepsilon =\const> 0$ and belonging to $\mathcal{C} (0, T;
\mathcal{D}' (\mathbb{R}_x^1))$ uniformly in $\varepsilon$, is
called a weak asymptotic mod $ O_{\mathcal{D}'}(\varepsilon^2)$
solution of (\ref{1}) if the relations
\begin{equation}\label{5}
\frac{d}{dt} \int_{-\infty}^\infty u\psi dx -
\int_{-\infty}^\infty  g'(u) \frac{\pa \psi}{\pa x} dx = O
(\varepsilon^2),
\end{equation}
\begin{equation}\label{6}
\frac{d}{dt} \int_{-\infty}^\infty u^2 \psi dx
+\int_{-\infty}^\infty \Big\{2\big(g(u)-ug'(u)\big)+
3 \left( \varepsilon\frac{\pa u}{\pa x}
\right)^2 \Big\}\frac{\pa \psi}{\pa x}dx = O (\varepsilon^2)
\end{equation}
hold uniformly in $t$ for any test function $\psi = \psi(x) \in
\mathcal{D} (\mathbb{R}^1)$.
\end{defi}
Here the right-hand sides are $\mathcal{C}^\infty$-functions for
$\varepsilon=\const > 0$ and  piecewise continuous functions
uniformly in $\varepsilon \geq 0$. The estimates are understood in
the $\mathcal{C}(0, T)$ sense:
$$g(t, \varepsilon) = O (\varepsilon^k) \leftrightarrow \max_{t\in[0,T]} |g(t,
\varepsilon)| \leq c\varepsilon^k.$$
\begin{defi} \label{def2.2}
A function  $v(t, x, \varepsilon)$ is said to be of the value $
O_{\mathcal{D}'}(\varepsilon^k)$ if the relation
$$ \int_{-\infty}^\infty v(t, x, \varepsilon )\psi(x) dx = O(\varepsilon^k) $$
holds uniformly in $t$ for any test function $\psi \in \mathcal{D}
(\mathbb{R}_x^1)$.
\end{defi}
Comparing the left-hand sides of (\ref{5}), (\ref{6}) with
(\ref{1}), (\ref{7}) we conclude that Definition 2 calls a
function the ``weak asymptotic solution" if it satisfies the
conservation laws (\ref{1}), (\ref{7}) in the sense $
O_{\mathcal{D}'}(\varepsilon^2)$.

Let us consider the interaction of two solitary waves for the
model (\ref{1}) with the initial data
\begin{equation} \label{8}
u|_{t=0} = \sum_{i=1}^2 A_i\omega \left(\beta_i \frac{x - x_i^0
}{\varepsilon},A_i\right),
\end{equation}
where
$A_2>A_1>0$, $x_1^0 -x_2^0=\const>0$
and we  assume  the same relations between $A_i$, $\beta_i$ and $V_i$ as in
(\ref{2h}).
Obviously, the trajectories $x = V_i t + x_i^0$
have a joint point $x = x^*$ at a time instant $t = t^*$.

Following \cite{DanOm, Om}, we write the asymptotic ansatz in  the  form:
\begin{equation}\label{9}
u = \sum_{i=1}^2 G_i(\tau)\omega \left( \beta_i \frac{x - \vp_i
(t, \tau, \varepsilon)}{\varepsilon},A_i  \right), \quad G_i(\tau)=A_i+S_i(\tau).
\end{equation}
Here $\vp_i = \vp_{i0}(t) + \varepsilon \vp_{i1}(\tau)$, where
$\vp_{i0} = V_i t + x_i^0$, are the trajectories of noninteracting
solitary waves;
$$\tau = \psi_0(t)/\varepsilon, \quad \psi_0(t) = \beta_1\big(\vp_{20}(t) - \vp_{10}(t)\big),$$
denotes the ``fast time"; the phase and amplitude corrections
$\vp_{i1}$, $S_{i}$ are smooth functions such that
\begin{equation}\label{10}
\vp_{i1}(\tau) \rightarrow 0 \quad \mbox{as} \; \tau \rightarrow
-\infty, \quad \vp_{i1} (\tau)\rightarrow \vp_{i1}^\infty =
\mbox{const}_i \quad \mbox{as} \; \tau \rightarrow +\infty,
\end{equation}
\begin{equation}\label{11}
S_i (\tau)\rightarrow 0 \quad \mbox{as} \; \tau \rightarrow \pm\infty
\end{equation}
with an exponential rate.
\subsection{Asymptotic construction}

To construct the asymptotics we should calculate  the weak expansions of the terms from the left-hand sides of the relations (\ref{5}), (\ref{6}).
It is easy to check that
\be
u=\ve \sum_{i=1}^2a_{1,i}\frac{G_i}{\beta_i}\delta(x-\vp_i)+O_{\mathcal{D}'}(\ve^3),\label{15}
\ee
where $\delta(x)$ is the Dirac delta-function. Here and in what follows we use the notation
\begin{equation}\label{15a}
a_{k,i}\ipd\int_{-\infty}^\infty\big(\om(\eta,A_i)\big)^k d\eta, \quad k>0,\qquad
a'_{2,i}\ipd\int_{-\infty}^\infty\big(\om'(\eta,A_i)\big)^2 d\eta.
\end{equation}
At the same time for any $F(u)\in C^1$
\begin{align}
&\int_{-\infty}^\infty F\left(\sum_{i=1}^2 G_i\omega \left( \beta_i \frac{x - \vp_i
}{\varepsilon},A_i  \right)\right)\psi(x)dx\notag\\&
=\ve\sum_{i=1}^2\frac{1}{\beta_i}\int_{-\infty}^\infty F\big(A_i\om(\eta,A_i)\big)
\psi(\vp_i+\ve\frac\eta\beta_i)d\eta+\frac{\ve}{\beta_2}\int_{-\infty}^\infty\Big\{ F\Big(G_1\om(\eta_{12},A_1)\label{16}\\
&+G_2\om(\eta,A_2)\Big)-F\big(A_1\om(\eta_{12},A_1)\big)-F\big(A_2\om(\eta,A_2)\big) \Big\}\psi(\vp_2+\ve\frac\eta\beta_2)d\eta,\notag
\end{align}
where
\be
\eta_{12}=\theta\eta-\sigma, \quad \sigma=\beta_1(\vp_1-\vp_2))/\ve\quad \theta=\beta_1/\beta_2.\label{17}
\ee
We take into account that the second integrand in right-hand side (\ref{16}) vanishes exponentially fast as $|\vp_1-\vp_2|$  grows,
thus, its main contribution is at the point $x^*$.  We write
\be
\vp_{i0}=x^*+V_i(t-t^*)=x^*+\ve \frac{V_i}{\dpsi}\tau \quad \text{and} \quad \vp_{i}=x^*+\ve\chi_{i},\label{18}
\ee
where $\dpsi=\beta_1(V_2-V_1)$, $\chi_i=V_i\tau/\dpsi + \vp_{i1}$.
It remains to apply the formula
\be
f(\tau)\delta(x-\vp_i)=f(\tau)\delta(x-x^*)-\ve\chi_if(\tau)\delta'(x-x^*)+O_{\mathcal{D}'}(\ve^2),\label{19}
\ee
which holds for each $\vp_i$ of the form (\ref{18}) with slowly increasing $\chi_i$ and for $f(\tau)$ from the Schwartz space.
 Moreover, the second term in (\ref{19}) is $O_{\mathcal{D}'}(\ve)$. Thus, under the assumptions (\ref{10}), (\ref{11})
  we obtain the weak asymptotic expansion of $F(u)$ in the final form:
\be
F(u)=\ve \sum_{i=1}^2F(A_i)\frac{a_{F,i}}{\beta_i}\delta(x-\vp_{i})+\ve\frac{F(A_2)}{\beta_2}\mathfrak{R}_{F}\delta(x-x^*)
+O_{\mathcal{D}'}(\varepsilon^2),\label{20}
\ee
where
\begin{align}
 a_{F,i}&=F(A_i)^{-1}\int_{-\infty}^\infty F\big(A_i\om(\eta,A_i)\big)d\eta, \notag\\
\mathfrak{R}_{F}&=
F(A_2)^{-1}\int_{-\infty}^\infty\Big\{ F\Big(G_1\om(\eta_{12},A_1)
+G_2\om(\eta,A_2)\Big)\label{21}\\
&-F\big(A_1\om(\eta_{12},A_1)\big)-F\big(A_2\om(\eta,A_2)\big) \Big\}d\eta.\notag
\end{align}
Note that to define $\pa u^2/\pa t \mod O_{\mathcal{D}'}(\varepsilon^2)$ it is necessary to
calculate $u^2$ with the precision $O_{\mathcal{D}'}(\varepsilon^3)$. Thus, transforming (\ref{15}) with the help
of (\ref{19}) and using (\ref{20}) with $F(u)=u^2$, we obtain modulo $O_{\mathcal{D}'}(\ve^3)$:
\begin{align}
& u=\ve \sum_{i=1}^2a_{1,i}K_{i0}^{(1)}\delta(x-\vp_i)
+\ve \sum_{i=1}^2a_{1,i}K_{i1}^{(1)}\Big\{\delta(x-x^*)-\ve\chi_i\delta'(x-x^*)\Big\},\notag\\
&u^2=\ve \sum_{i=1}^2a_{2,i}K_{i0}^{(2)}\delta(x-\vp_i)
+\ve \sum_{i=1}^2\Big\{a_{2,i}K_{i1}^{(2)}+2\tilde{a}_2\frac{G_1G_2}{\beta_2}R_2^{(0)}\Big\}\delta(x-x^*)\notag\\
&-\ve^2\Big\{\sum_{i=1}^2a_{2,i}K_{i1}^{(2)}\chi_i+2\tilde{a}_2\frac{G_1G_2}{\beta_2}\Big(\chi_2R_2^{(0)}
+\frac{1}{\beta_2}R_2^{(1)}\Big)\Big\}\delta'(x-x^*),\label{15b}
\end{align}
where
\begin{align}
&
K_{i}^{(n)}=\frac{G^n_i}{\beta_i},\quad K_{i0}^{(n)}=\frac{A^n_i}{\beta_i},
\quad K_{i1}^{(n)}=K_{i}^{(n)}-K_{i0}^{(n)},\label{15c}\\
&R_2^{(i)}=\frac{1}{\tilde{a}_2}\int_{-\infty}^\infty\eta^i\om(\eta_{12},A_1)\om(\eta,A_2)d\eta,
\quad \tilde{a}_2=\sqrt{a_{2,1}a_{2,2}}.\label{15d}
\end{align}
Calculating weak expansions for other terms from the Definition 2  and substituting them into
(\ref{5}), (\ref{6}) we obtain linear combinations of
$\delta'(x-\vp_{i})$, $i=1,2$, $\delta(x-x^*)$, and
$\delta'(x-x^*)$ (see also \cite{DanOm, DanOmShel}). Therefore, we obtain:

system of algebraic equations
\begin{align}
&a_{1,i}V_iK_{i0}^{(1)}-a_{g',i}\,g'(A_i)/\beta_i=0, \quad i=1,2,\label{22}\\
&a_{2,i}V_iK_{i0}^{(2)}+2a_{g_2,i}\,g_2(A_i)/\beta_i+3a'_{2,i}\beta_i^2K_{i0}^{(2)}=0, \quad i=1,2,\label{23}
\end{align}
system of functional equations
\be
\sum_{i=1}^2a_{1,i}K_{i1}^{(1)}=0,\quad \sum_{i=1}^2a_{2,i}K_{i1}^{(2)}+2\tilde{a}_2\frac{G_1G_2}{\beta_2}{R}^{(0)}_2=0,\label{24}\\
\ee
and system of ordinary differential equations
\begin{align}
\dpsi\frac{d}{d\tau}\sum_{i=1}^2a_{1,i}\Big\{K_{i0}^{(1)}\vp_{i1}&+\chi_iK_{i1}^{(1)}\Big\}=f,\label{25}\\
\dpsi\frac{d}{d\tau}\Big\{\sum_{i=1}^2a_{2,i}\Big(K_{i0}^{(2)}\vp_{i1}&+\chi_iK_{i1}^{(2)}\Big)\label{26}\\
&+2\frac{\tilde{a}_{2}}{\beta_2}G_1G_2\Big(\chi_2{R}^{(0)}_2+\frac{1}{\beta_2}{R}^{(1)}_2\Big)\Big\}=F,\notag
\end{align}
where
\begin{align}
&f=\frac{g'(A_2)}{\beta_2}\mathfrak{R}_{g'},\quad g_2=g(u)-ug'(u),\quad \tilde{a}'_{2}=\sqrt{a'_{2,1}a'_{2,2}},\label{35}\\
& F=-2\frac{g_2(A_2)}{\beta_2}\mathfrak{R}_{g_2}-3\left\{\sum_{i=1}^2a'_{2,i}\beta_i^2K_{i1}^{(2)}
+2\tilde{a}'_{2}\beta_1G_1G_2{R}^{(0)}_{2,1}\right\},\label{35a}
\end{align}
 and ${R}^{(0)}_{2,1}$ is of the form (\ref{15d})
for ${R}^{(0)}_{2}$ but with $\tilde{a}'_{2}$, $\om'$ instead of $\tilde{a}_{2}$, $\om$.
\begin{lem}
The algebraic equations (\ref{22}), (\ref{23}) imply again the  relations
(\ref{2h}) between $A_i$, $\beta_i$, and $V_i$.
\end{lem}
{\bf Proof}.  Let us change $A\to A_i$, $V\to V_i$, $\beta\to\beta_i$ in (\ref{2f}) and
integrate it with respect to $\eta$.   Then we obtain the equality
 \be
 a_{1,i}A_iV_i=a_{g',i}\,g'(A_i),\label{27}
 \ee
which is equivalent to (\ref{22}). Next let us multiply the original third order ordinary equation for $\om$ by $\om$.
Integrating we obtain the following alternative version of (\ref{2f}):
$$
V\om^2+2\frac{g_2(A\om)}{A^2}=\beta^2\Big\{\frac{d^2(\om^2)}{d\eta^2}-3\left(\frac{d\om}{d\eta}\right)^2\Big\}.
$$
Thus
\be
 a_{2,i}V_i+2a_{g_{2},i}\frac{g_{2}(A_i)}{A_i^2}=-3 a'_{2,i}\beta_i^2.\label{28}
 \ee
At the same time, integration of (\ref{2g}) implies:
\be
 a_{2,i}V_i-2a_{g,i}\frac{g(A_i)}{A_i^2}=a'_{2,i}\beta_i^2.\label{28a}
 \ee
Taking into account (\ref{28}) we obtain the relation
\be
 a_{g,i}\,g(A_i)=-a_{g_{2},i}\,g_2(A_i)-2a'_{2,i}\beta_i^2A_i^2,\label{28b}
 \ee
which implies that (\ref{23}) is the result of the integration of (\ref{2g}) with respect to $\eta$.
Now we square the both parts of (\ref{3}) and integrate them:
\be
 a'_{2,i}=a_{2,i}-a_{g,i}.\label{28c}
 \ee
Substituting this into (\ref{28a}) we obtain the equality
$$
a_{2,i}\big(V_i-\beta_i^2\big)-a_{g,i}\,g_1(A_i)\Big(2-\frac{\beta_i^2}{g_1(A_i)}\Big)=0,
$$
which implies the relations (\ref{2h}).
\qed
\medskip
\subsection{Analysis of the model equations  (\ref{24}) - (\ref{26})}
It is easy to note that the  system (\ref{24}) implies the quadratic equation with
coefficients which depend on the convolution $R^{(0)}_2=R^{(0)}_2(\sigma)$.   To analyze the  equation we assume:
\be\label{29a}
A_i>>1,\quad i=1,2,\qquad\theta<<1.
\ee
Moreover, let
\be\label{29b}
g_1(z)=\sum_{k=1}^nc_kz^{q_k},\quad \delta_1\leq q_1<q_2<\dots<q_n<4.
\ee
Then (\ref{2h}) imply the relation
\be\label{29c}
A_i=c'\beta_i^{q'}\Big(1+O\big(A_i^{q_{n-1}-q_n}\big)\Big),\quad q'=2/q_n,\; c'=(2c_n)^{-1/q_n},\; i=1,2.
\ee
Let us define the following notation:
\be
\kappa_i=\frac{S_i}{\beta_i}\Big(c'\beta_2^{q'-1}\Big)^{-1}, \quad i=1,2.\label{29}
\ee
\begin{defi} \label{def2.3}
A function  $f(\tau,\theta)$ is said to be of the value $O_S(\theta^k)$ if there exists  a function
$s(\tau)$ from the Schwartz space such that the estimate
\be
|f(\tau,\theta)|\leq c\,\theta^k|s(\tau)| \label{30}
\ee
holds uniformly in $\tau$ for a constant $c>0$.
\end{defi}
Let us note also that, under the condition (\ref{10}), the convolution $R^{(0)}_2\to0$ as $\tau\to\pm\infty$ with an exponential rate.

Assumptions (\ref{29a}), (\ref{29b}) allow to prove the statement:
\begin{lem}
Let the assumptions (\ref{29a}), (\ref{29b}) be satisfied.  Then the  algebraic equations (\ref{22}), (\ref{23})
allow to define $S_i$ with the
property (\ref{11}). More in detail,
 for sufficiently small $\theta$ we find
\be
\kappa_1=\frac{\sqrt{\ol{a}_2}}{\ol{a}_1}\theta^{q'}R^{(0)}_2+O_S(\theta^2+\theta^{2q'}),\quad \kappa_2=-\ol{a}_1\kappa_1,
\quad\ol{a}_i=a_{i,1}/a_{i,2}.\label{31}
\ee
\end{lem}
 Now let us simplify the equations (\ref{25}), (\ref{26}). We note firstly that
in view of the first equation (\ref{24}) and the identity
\be
\beta_1(\chi_1-\chi_2)=\sigma \label{32}
\ee
it is possible to eliminate $\chi_i$ from the left-hand side of (\ref{25}), since
$$
\sum_{i=1}^2a_{1,i}\big\{K_{i0}^{(1)}\vp_{i1}+\chi_iK_{i1}^{(1)}\big\}=\sum_{i=1}^2a_{1,i}K_{i0}^{(1)}\vp_{i1}+a_{1,1}\frac{\sigma}{\beta_1}K_{11}^{(1)}.
$$
Simplifying in the same manner   the equation (\ref{26}),
 we transform (\ref{25}), (\ref{26}) to the following form:
\begin{align}
&\dpsi\frac{d}{d\tau}\Big\{\sum_{i=1}^2a_{1,i}K_{i0}^{(1)}\vp_{i1}+a_{1,1}\frac{\sigma}{\beta_1}K_{11}^{(1)}\Big\}=f,\label{33}\\
&\dpsi\frac{d}{d\tau}\Big\{\sum_{i=1}^2a_{2,i}K_{i0}^{(2)}\vp_{i1}+a_{2,1}\frac{\sigma}{\beta_1}K_{11}^{(2)}
+2\tilde{a}_{2}\theta K_{1}^{(1)}K_{2}^{(1)}{R}^{(1)}_2\Big\}=F,\label{34}
\end{align}
where $f$ and $F$ are defined in (\ref{35}), (\ref{35a}).

The second step is the elimination of $\vp_{i1}$ from the model system. To do it we divide $\sigma$ into
the growing ${\beta_1}(V_1-V_2)\tau/{\dpsi}=-\tau$ and the bounded
(if the assumptions (\ref{10}) are satisfied) $\ts=\sigma+\tau$ parts. Since
\be\label{36}
\vp_{11}=\vp_{21}+\tilde{\sigma}/\beta_1,
\ee
we obtain from (\ref{33})
\be
\dpsi\frac{d}{d\tau}\left\{{r_1}\vp_{21}+\frac{\sigma}{\beta_1}K_{1}^{(1)}
\right\}=\frac {f}{a_{1,1}}-\frac{\dpsi}{\beta_1}K_{10}^{(1)}.\label{37}
\ee
Here and in what follows we use the notation
\be
r_j=\sum_{i=1}^2\frac{a_{j,i}}{a_{j,1}}K_{i0}^{(j)}\quad \text{for}\quad  j=1\quad \text{and}\quad j=2.\label{38}
\ee
Now,  transforming  (\ref{34}) in the same manner  and applying the first assumption (\ref{10}) we pass to the problem:
\be
\frac{d}{d\,\tau}\, Q(\sigma)=\mathfrak{F}(\sigma),\quad \frac{\sigma}{\tau}\Big|_{\tau\to-\infty}\to-1, \label{39}
\ee
where
\begin{align}
&Q=\frac{\sigma}{\beta_1}\Big\{K_{1}^{(2)}-\frac {r_2}{r_1}K_{1}^{(1)}\Big\}+\frac{2}{\sqrt{\ol a}_2}K_1^{(1)}K_2^{(1)}R_2^{(1)},\notag\\
&\mathfrak{F}=-\frac{1}{\beta_1}\Big\{K_{10}^{(2)}-\frac {r_2}{r_1}K_{10}^{(1)}\Big\}+\frac{1}{\dpsi}\Big\{\frac{1}{a_{2,1}}F-\frac {r_2}{a_{1,1}r_1}f\Big\}.\notag
\end{align}
Sufficiently simple analysis of the equation  (\ref{39}) implies the  statement:
\begin{lem}
Under the assumptions (\ref{29a}), (\ref{29b}) the following relations hold:
\begin{align}
&\frac{dQ}{d\,\sigma}=-\frac{A_1A_2}{\beta_1^2}\Big\{\frac{\bar{a}_2}{\bar{a}_1}-\theta^{q'}+O_S(\theta+\theta^{q'})\Big\},\label{40} \\
&\mathfrak{F}=\frac{A_1A_2}{\beta_1^2}\Big\{\frac{\bar{a}_2}{\bar{a}_1}-\theta^{q'}+O_S(\theta+\theta^{q'})\Big\}.\label{40a}
\end{align}
\end{lem}
The uniform in $\tau$ inequality $\mathfrak{F}>0$ and the exponential type behavior of $\mathfrak{F}$
and $Q$ imply the existence of the function $\sigma$ such that $\ts=\sigma+\tau$ is bounded and tends to
its limiting values with an exponential rate. This and the equalities (\ref{36}), (\ref{37})
justify the existence of the required phase corrections $\vp_{i1}$
with the property (\ref{10}).
\qed
\medskip

Our  main result is the following:
\begin{teo} \label{teo1}
Let the assumptions (\ref{29a}), (\ref{29b}) be satisfied.
Then the solitary wave collision  in the problem
(\ref{1}), (\ref{8}) preserves the elastic scenario with accuracy
$O_{\mathcal{D}'} (\varepsilon^2)$ in the sense of Definition
\ref{def2.1}.  The weak asymptotic solution  has the form
(\ref{9}).
\end{teo}
The next theorem allows to treat the weak asymptotics (\ref{22}) in the classical sense:
\begin{teo} \label{teo00}
Let the assumptions (\ref{29a}), (\ref{29b}) be satisfied. Then the function $u$ of the form (\ref{9})
is a weak asymptotic mod $ O_{\mathcal{D}'}(\varepsilon^2)$ solution of (\ref{1})
if and only if $u$ satisfies the following conservation and balance laws:
\begin{align}
&\frac{d}{dt}\int_{-\infty}^\infty u\,dx=0,\quad \frac{d}{dt}\int_{-\infty}^\infty u^2dx=0,\label{12}\\
&\frac{d}{dt}\int_{-\infty}^\infty xu\, dx-\int_{-\infty}^\infty
g'(u) dx=0,\label{13}\\
&\frac{d}{dt}\int_{-\infty}^\infty xu^2
dx+2\int_{-\infty}^\infty g_2(u)dx+3
\int_{-\infty}^\infty\left(\ve\frac{\pa u}{\pa x}\right)^2dx=0.\label{14}
\end{align}
\end{teo}
To prove the Theorem 3 it is enough to rewrite the equalities
  (\ref{22}) - (\ref{26}) as integrals of the function
 (\ref{9}) and its derivatives. Results of numerical simulations confirm the traced asymptotic
analysis (see \cite{OmVal} for the nonlinearity $u^{3/2}$).
\section{Dynamics of perturbed solitary waves}

In this section  we  consider  briefly the perturbed KdV-type equation (\ref{1}),
\begin{equation}\label{41}
\frac{\pa u}{\pa t}+\frac{\pa g'(u)}{\pa x}+\varepsilon^2
\frac{\partial^3u}{\pa x^3}= F,
\end{equation}
where $F=F(x,t,u,\ve u_x, \ve^2 u_{xx},\dots)\in C^\infty$ is  ``small" in our scaling.
We assume that $F|_{u\equiv0}=0$.

Let us construct firstly a self-similar one-phase asymptotic solution   and discuss after that how
to use this asymptotics for more realistic Cauchy data.

According to results \cite{Om, Om2} to construct the leading term of the classical one-phase asymptotic solution it is enough to  find the weak
asymptotics. By analogy  with  Definition 2 we write:
\begin{defi}
 A sequence $u(t, x, \varepsilon )$, belonging to
$\mathcal{C}^\infty (0, T; \mathcal{C}^\infty (\mathbb{R}_x^1))$
for $\varepsilon =\const> 0$ and belonging to $\mathcal{C} (0, T;
\mathcal{D}' (\mathbb{R}_x^1))$ uniformly in $\varepsilon$, is
called a weak asymptotic mod $ O_{\mathcal{D}'}(\varepsilon^2)$
solution of (\ref{41}) if the relations
\begin{align}
&
\frac{d}{dt} \int_{-\infty}^\infty u\psi dx -
\int_{-\infty}^\infty  g'(u) \frac{\pa \psi}{\pa x} dx -\int_{-\infty}^\infty F\psi dx= O(\varepsilon^2),\label{41a}\\
&
\frac{d}{dt} \int_{-\infty}^\infty u^2 \psi dx
+\int_{-\infty}^\infty \Big\{2g_2(u)+
3 \left( \varepsilon\frac{\pa u}{\pa x}
\right)^2 \Big\}\frac{\pa \psi}{\pa x}dx\label{41b}\\
&\qquad\qquad\qquad\qquad\qquad\qquad\quad -2\int_{-\infty}^\infty Fu\psi dx= O (\varepsilon^2)\notag
\end{align}
hold uniformly in $t$ for any test function $\psi = \psi(x) \in
\mathcal{D} (\mathbb{R}^1)$.
\end{defi}
Combining the ideas of  \cite{MasOm1, MasOm2} and \cite{Om, Om2} we write the ansatz in the form:
\begin{equation}\label{42}
u=A\omega \big( \beta(x-\vp(t))/\ve, A \big)+\ve Y(\tau,t,x),
\end{equation}
where $A=A(t)$, $\beta=\beta(t)$,
 and $Y$ is a smooth bounded function such that
 $Y(\tau,t,x)\to0$ as $\tau\to+\infty$,
and $Y(\tau,t,x)\to u^-(x,t)$ as $\tau\to-\infty$. Note that (\ref{42}) can be treated as a "two-phase"
asymptotics since
\be\label{42b}
u=\ve a_{1}\frac{A}{\beta}\delta\big(x-\vp(t)\big)+\ve u^-(x,t)H\big(\vp(t)-x\big)+O_{\mathcal{D}'}(\ve^2),
\ee
 and the coefficient of the Heaviside function $H$ varies slowly.

 We take into account the relation
\be
F=\ve \frac{a_{F_0}}{\beta}{\bar F}\delta\big(x-\vp(t)\big)
+\ve F'_u|_{u\equiv0}u^-(x,t)H\big(\vp(t)-x\big)+O_{\mathcal{D}'}(\ve^2),\label{42c}
\ee
where ${\bar F}=F(\vp,t,A,\beta A,\dots)$, $F_0=F(\vp,t,A\om,A\beta\om', A\beta^2\om'',\dots)$. Next calculating others weak expansions   and substituting them into
(\ref{41a}), (\ref{41b}) we obtain linear combinations of
$\delta(x-\vp)$,  $\delta'(x-\vp)$, and
$H(\vp-x)$. Therefore, we pass to the following system:
\begin{align}
&a_1A\frac{d\vp}{dt}=a_{g'}g'(A),\quad a_2\frac{d\vp}{dt}+2a_{g_2}\frac{g_2(A)}{A^2}+3a_2'\beta^2=0,\quad t>0,\label{43}\\
&\frac{d}{dt}\left(a_2\frac{A^2}{\beta}\right)=2a_{\om F_0}\frac{A}{\beta}{\bar F},\quad t>0,\label{44}\\
&\frac{\pa u^-(x,t)}{\pa t}=F'_u(x,t,0,\dots)u^-(x,t),\quad t>0,\quad x<\vp(t),\label{45}\\
&u^-(\vp,t)\frac{d\vp}{dt}+\frac{d}{dt}\left(a_1\frac{A}{\beta}\right)=\frac{a_{F_0}}{\beta}{\bar F}\quad t>0.\label{46}
\end{align}
 Lemma 1 implies that the  equations (\ref{43}) are equivalent to the equalities
\be
\frac{d\vp}{dt}=\beta^2,\quad \beta^2=2g_1(A).\label{47}
\ee
Thus, the equations   (\ref{44}), (\ref{47}) form the complete system to define $\vp$, $A$, and $\beta$.

Next note that    $u^-$ is defined by (\ref{45}) for $x<\vp(t)$ whereas
 the equation (\ref{46}) defines the boundary value of $u^-$ on the curve $x=\vp(t)$.

Let us stress finally that the self-similarity implies a special choice of the initial data. In particular,
the initial function $Y(\tau,0,x)$
should be of the special form
\begin{equation}\label{47b}
Y(\tau,0,x)=\Big\{u_0^-(x)\chi(\tau,t)+Z_1(\tau,t)+c_1\om'(\tau,A)\Big\}\Big|_{t=0},
\end{equation}
where $u_0^-$ is a smooth function such that $u_0^-\big(\vp(0)\big)=u^-\big(\vp(0),0\big)$,
$\chi(\tau,t)$ is a  regularization of the Heaviside function,  $Z_1(\tau,t)$ is a special function from the Schwartz space,
and  $c_1$ is  arbitrary constant (see \cite{MasOm1, MasOm2}).
If it is violated and, for example,
\be u|_{t=0}=A(0)\om\Big(\beta \big(x-\vp(0)\big)/\varepsilon,A(0)\Big),\label{47c}
\ee
then the perturbed soliton generates a rapidly oscillating tail of the
amplitude $o(1)$ (the so called ``radiation") instead of the smooth tail $\ve u^-(x,t)$
(see \cite{Kal} for the perturbed KdV equation and numerical results in \cite{GO, GO2}). However,
$\ve u^-(x,t)$ describes sufficiently well the tendency of the radiation amplitude behavior.

\textbf{Example 1} \cite{OmVal}. Let $g'(u)=u^{3/2}$ and let
\be
F=-\frac\ve{2b}\frac{\pa}{\pa x}u^2+\ve\frac{b}{2}\frac{\pa}{\pa t}\Big\{\frac{u^2}{2b}-u^{3/2}-\ve^2\frac{\pa^2 u}{\pa x^2}
-\frac12\frac{\pa}{\pa t}\int_{-\infty}^xu\,dx'\Big\}, \; b=\const. \label{47a}
\ee
This right-hand side represents the remainder which was omitted in \cite{Sch} in the process of the regular
asymptotic construction. At the same time, for large $x$ and $t$ singular perturbations can appear and we should estimate
the influence of (\ref{47a}) on the solitary wave. However it is easy to check that $a_{\om F_0}=0$, thus $A=\const$.
Consequently, $u^-(\vp,t)=0$ and $u^-(x,t)\equiv0$. This justifies the elimination of $F$ from the leading terms of the asymptotics.

\textbf{Example 2} \cite{OmVal}. Let
\be
g'(u)=u^{3/2},\quad F=\mu(\alpha-u)u, \label{48}
\ee
where $F$ describes an external force, $\mu>0$, $\alpha>0$ are constants. Then the equation  (\ref{44}), supplied by the initial condition, takes the form:
\be
\frac{d A}{dt}=\mu' A\Big(1-\frac{A}{A^*_\alpha}\Big),
\quad A|_{t=0}=A^0,\label{49}
\ee
where $\mu'=8\alpha\mu/7$ and $A^*_\alpha=\alpha\, a_2/a_3$. Integrating we obtain
\be
A=A^0e^{\mu' t}\{1+c_\alpha(e^{\mu' t}-1)\}^{-1}, \quad c_\alpha=A^0/A^*_\alpha.\label{50}
\ee
Therefore, when $t$ grows, $A$ decreases to $A^*_\alpha$  if $A^0\geq A^*_\alpha$ or $A$ increases to $A^*_\alpha$  if $A^0\leq A^*_\alpha$.
Next integrating the equation (\ref{47}),
we conclude that the curve $x=\vp(t)$ tends to a straight line as $t\to\infty$. Thus, the solitary wave demonstrates a stable behavior.

Let us turn to the correction $u^-$. Preserving the term $O(\ve)$, we write the equation (\ref{45}) as following:
\be
\frac{\pa u^-}{\pa t}=\mu u^-(\alpha-\ve u^-).\label{51}
\ee
Simple algebra justifies that  $u^-(\vp,t)=O(\mu)>0$. Thus,   the amplitude of $u^-$ increases exponentially fast with the rate $O(\mu)$.
Moreover, it tends to the value $O(1/\ve)$,
so that the correction $\ve u^-$ becomes of the same value as the leading term  in a critical time $T^*\sim\ln(1/(\ve\mu))/(\alpha\mu)$.

Results of numerical simulation confirm this analysis. Namely, we consider the Cauchy problem  (\ref{41}), (\ref{47c}), (\ref{48}).
Since the initial value does not include the correction of the form (\ref{47b}),
 the soliton correction is not a smooth tail $\ve u^-(x,t)$,  but the radiation. However,
the behavior of the correction's amplitude  is explicitly the same as it been described above, see \cite{OmVal}.
\section{Conclusion}
In fact, the result that each equation from the family (\ref{1}) preserves the KdV-type scenario of soliton
interaction was rather expected. It was found much more interesting to consider the behavior of perturbed solitary waves.
It  turned out that there exists a class of perturbations  which provoke a very interesting scenario of instability development:
a short solitary wave (with the wave-length $\sim\ve$) varies its parameters, remaining stable but
generating a long wave  (with the wave-length $\sim\ve^\nu$, $\nu<1$) perturbation of a small amplitude. Inversely,
this perturbation turns out to be unstable, its amplitude increases and destroys the original soliton.

On the other hand, the rate of the perturbation growth can be slow (of the order $O(\mu)$ with $\mu<<1$ for the external
force (\ref{48})), whereas the amplitudes of solitons tend to the same stationary value $A^*_\alpha$.
So, for sufficiently small $\mu$ (for sufficiently large distances between the original positions of the solitons),
the amplitudes can be almost of the same value before the collision of solitons, which prevents the intersection of the trajectories.
In other words, the perturbation can eliminate the interaction between solitons. To illustrate the situation we refer to \cite{OmVal},
where the  dynamics of five  solitons is depicted for the equation (\ref{41}), (\ref{48}) with $\mu=0.2$ and the following amplitudes of the
original solitons:  $A_5=4$, $A_4=2$, $A_3=1$, $A_2=0.5$, and $A_1=0.25$.

\section{Acknowledgement}
The research was supported by  SEP-CONACYT under grant~178690 (Mexico).

\end{document}